\documentclass[12pt]{article}

\usepackage{epsf,epsfig,amsfonts,graphicx,a4wide} 
\usepackage{url}
\usepackage{amssymb}
\usepackage{amsthm}
\usepackage{epsf,amsfonts,graphics}
\usepackage{graphicx} 
\usepackage{a4wide,epsfig}
\newtheorem{theorem}{Theorem}
\newtheorem*{lemma}{Lemma}

\newcommand{\weg}[1]{}

\title{Locally 2-fold symmetric manifolds are locally symmetric} 
\author{Shaoqiang Deng and  Vladimir S. Matveev\thanks{ The result was obtain during the visit of V.M. to Nankai University related to the 2016 International Conference on Riemann-Finsler Geometry; he thanks the Nankai University and the conference organizers for their hospitality and financial support; V.M. also acknowledges the financial support of DFG. We are grateful to O. Yakimova for useful discussions. }}  
\date{} 
\begin{document} 

\maketitle

\abstract{
A manifold is locally \emph{$k$-fold symmetric}, if  for any point and any $k$-dimensional  vector   subspace  tangent to this point there exists a local  
 isometry such that this point is a fixed point and the differential of the isometry restricted to that   $k$-dimensional  vector   subspace is minus the identity. We show that for $k \ge 2$, Riemannian, pseudoriemannian and Finslerian  locally $k$-fold symmetric manifolds are locally symmetric. } 
\section*{Introduction} 

A Riemannian or Finslerian manifold $(M, g)$ or $(M,F)$ will be called  locally \emph{$k$-fold symmetric}, if for any point $x$ and for any $k$-dimensional vector subspace  $V^k\subseteq T_xM$   there exists $r>0$ and 
 a local isometry $\phi:B_r(x)\to B_r(x)$ of the $r-$ball around $x$ such that $f(x)= x$ and  its differential at the point $x$ restricted to $V^k$ is minus the identity.   The definition survives in the pseudoriemannian setting (when the metric $g$ is nondegenerate but  may have  arbitrary signature) if we merely allow the isometry to be defined in some neighborhood of $x$ and not necessary in the $r$-ball. 
 
If $k= \textrm{dim}(M),$ then the definition above is one of the standard definitions of locally-symmetric spaces; locally $1$-fold symmetric spaces are called \emph{ locally weakly symmetric}.

Our goal it to prove the following two statements.

\begin{theorem} \label{thm1}
For $k\ge 2$, each  locally $k$-fold symmetric Riemannian or pseudoriemannian manifold   is locally symmetric.  
\end{theorem}

\begin{theorem} \label{thm2}
For $k\ge 2$, each locally  $k$-fold symmetric Finslerian  manifold is locally symmetric and  Berwaldian. 
\end{theorem} 

Since  each  locally $k+1$-fold   symmetric spaces is clearly also  a locally $k$-fold symmetric, it is sufficient to prove both statements  for $k=2$ only.   
For $k=1$ the  statements  are wrong, see \cite{Deng} for a counterexample. In the global setting (i.e., when the manifold is complete), both statements were proved in 
\cite{Deng}; the proof is quite involving, is essentially global,  and in particular extensively uses the quite nontrivial results  on weakly symmetric  spaces from \cite{Yakimova0,Yakimova1,Yakimova2,Yakimova3}. Our proof is much easier,  shorter, and works locally.   

Let us explicitly mention that in the definition of locally $k$-fold symmetric metrics we do not require that the radius $r= r(x, V^k)$ from the definition  is locally bounded from below; in fact, if it is locally bounded from below then  (in the Riemannian and Finlser cases)  the algebraic methods of \cite{Deng} and \cite{Yakimova0,Yakimova1,Yakimova2,Yakimova3} are still applicable. But a consequence of our results is that locally $(k\ge 2)$-fold symmetric spaces are locally isometric  to globally symmetric spaces so that the radius $r= r(x, V^k)$ is locally bounded from below.  

All object in our paper are assumed to be sufficiently smooth; $C^3$-smoothness is enough. The $k$-fold symmetric Finlser  metrics  $F$ are automatically reversible.    Theorem \ref{thm2} remains true in the case of not necessary strictly  convex Finsler  metrics: indeed, similar to the end of the proof of \cite[Theorem 9.2]{MT2012} one can  show that if the theorem is true for strictly convex Finsler metrics it is true for arbitrary  Finsler metrics.

\section{Proof of Theorem \ref{thm1}}

We denote by $\nabla R$ the covariant derivative of the Riemann curvature tensor with lowered indexes, $\nabla R$ is a $(0,5)$-tensorfield $$\nabla R(U,V,X,Y,Z)= R_{ijk\ell,m} U^iV^jX^kY^\ell Z^m.$$ Our goal is to show that $\nabla R=0$, since Riemannian and pseudoriemannian metrics such that $\nabla R=0$ are locally symmetric. 

 By the  definition of locally $2$-fold symmetric spaces, for any  $x\in M$ and any two vectors $U, V \in T_xM$ there exists an isometry  
$\phi$ of $B_r(x)$ such that   $\phi(x)= x$ and  \begin{equation} \label{eq1} d\phi_x(U)= -U\ , \  d\phi_x(V)= -V.\end{equation} We use this isometry and  apply the Cartan's trick: 
since  $\nabla R$ is preserved by the isometries, (\ref{eq1}) implies   
$$\begin{array}{r} \nabla R(U,V,U, V,  U)= \nabla  R \left( d\phi_x(U),d\phi_x(V),d\phi_x(U), d\phi_x(V),  d\phi_x(U)  \right)\\ \stackrel{(\ref{eq1})}{=}\nabla  R(-U, -V,  -U, - V, -U)= - \nabla R(U,V,U, V,  U), \end{array} $$ which implies that for any two vectors $U,V$ we have    $\nabla R(U,V,U, V,  U)= 0$. 

Replacing  $U$ by $U+ t X + s Y $ and $V$ by $V+ r Z$ in the left-hand side of the last equation we obtain a  polynomial function in $t, s, r$; vanishing the coefficient of $ts r$ implies in view of the algebraic symmetries of the covariant derivative of the   Riemann curvature tensor the equation 
\begin{equation} \label{eq} \begin{array}{l} 
\nabla R(X, V, Y, Z, U) + \nabla R(Y, V, U, Z, X) +  \nabla R(U, V, X, Z, Y) + \\ \nabla R(X, Z, Y, V, U)   + \nabla R(Y, Z, U, V, X) +  \nabla R(U, Z, X, V, Y)=0.\end{array}   
\end{equation} 

Since vanishing of $\nabla R$ implies that the metric is locally symmetric, 
Theorem \ref{thm1} follows from the following linear algebraic lemma: 

\begin{lemma} \label{lemma}
Suppose a $(0,5)$-tensor $\nabla R=R_{ijk\ell , m}$ has the algebraic symmetries of the covariant derivative of the Riemann curvature tensor (i.e., is skewsymmetric  with respect to first two indexes, symmetric with respect to first and second pairs of indexes, satisfies the Bianchi equality with respect to the first three and with respect to the last three indexes). If in addition (\ref{eq}) is fulfilled, then it vanishes. 
\end{lemma}

{\bf Proof. } The  equation (\ref{eq}) is a system of linear equations on the components of $\nabla R$ and its solution space is the kernel  of  the linear mapping 
\begin{equation} \label{eqmap}
  \nabla R  \mapsto T, \ T_{ijk\ell m}:=  R_{ijk\ell, m} + R_{kjm\ell, i} + R_{mji\ell, k} + R_{i\ell kj , m}  + R_{k\ell m j , i} + R_{m\ell i j , k}. 
\end{equation}
We show that this mapping has trivial kernel, or equivalently, that  the  mapping  (\ref{eqmap}) has a left-inverse mapping. By direct calculation,  we see that the mapping 
 \begin{equation} \label{eqmap1} T\mapsto \nabla  R,  \ R_{{j i m k,  \ell }}=-\frac{1}{6} T_{i j \ell  k m}-\frac{1}{12} { T}_{{j k i \ell
m}}+\frac{1}{12}\,{ T}_{{j \ell  i m k}}+\frac{1}{6} { T}_{{i j k m \ell}}.
  \end{equation} 
 is left inverse to  (\ref{eqmap}), that is  the superposition $$  \nabla R \stackrel{(\ref{eqmap})}{\mapsto} T  \stackrel{(\ref{eqmap1})}{\mapsto} \nabla R       $$ is the identity.    Lemma and therefore Theorem \ref{thm1} are  proved. 

\section{Proof of Theorem  \ref{thm2}}

We assume that  our Finsler metric $F$ is locally 2-fold symmetric; our first goal will be  to show that  it is \emph{Berwaldian}, that is, there exists a torsion-free  affine connection $\nabla =\left(\Gamma^i_{jk}\right)$ (called  \emph{the associated connection}) such that the parallel transport in this connection preserves the Finlser metric.  

The  proof is  based on the same ideas as the   proof  of \cite[Theorem 9.2]{MT2012},  and will use the Binet-Legendre metric $g_F$ introduced in 
\cite[Definition 1.1]{MT2012}, we recommend to have this paper by hand.   In fact we are going to show that the Levi-Civita connection of $g_F$ is the associated connection of $F$. 
It is sufficient to show that the spray coefficients $G^i(x, V)$ of $F$ 
and the Christoffel symbols  of $g_F$ are related by  
$$
G^i(V)=  \frac{1}{2} \Gamma^i_{jk} V^jV^k.  
$$

In order to do it, we consider  an arbitrary point $x_0 \in M$ and work in the $g_F$-exponential  coordinate system around $x_0$; in this coordinate system straight lines passing through $x_0$ are $g_F
$-geodesics and  $\Gamma(x_0)^i_{jk}\equiv 0$; our goal is to show that for any $V\in T_{x_0}M$ we have $G^i(V)=0$. Consider the  Finslerian geodesic 
$\gamma(t)$ such that $\gamma(0)=x_0$ and $\dot\gamma(0)= V$. It is known  that  every  $F$-geodesic satisfies the  equation 
$
\ddot \gamma^i (t) + 2 G^i(\gamma(t),\dot\gamma(t))=0$, therefore at $x_0$ we have  
$$
\ddot \gamma^i (0) + 2 G^i(x_0,V)=0.
$$

Would $G^i(x_0,V)\ne 0$, at least one of the second derivatives of $\gamma^i(0)$ will be not zero and  therefore 
the geodesic $\gamma$ would  have at most two  intersections with any  straight line passing through $x_0$.
Let  us show that it has at least three intersections with some straight lines.

In order to do it, for $U\in T_xM$ we denote   by $\phi^{U}$ a local 
 isometry from the definition of locally  $k$-fold symmetric spaces such that $\phi^U(x_0)= x_0$, $d\phi_{x_0}^U(\dot\gamma(0))=  -\dot \gamma(0)$ 
and $d\phi^U_{x_0}(U)= -U$. Since the sphere is compact, 
 there exists a finite number  $\phi^{U_1},...,\phi^{U_N}$ of  such isometries  such that for any $Z\in T_{x_0}M$ there exists $i\in 
\{1,...,N\}$   such that $d\phi^{U_i}_{x_0}(Z)\ne Z$. We can think  that all these isometries operate on the ball of radius $r>0$ around $x_0$.

Take very small  $\varepsilon>0$ which in particular is much less  
 than $r$  and consider the points $\gamma(-\varepsilon)$ and $\gamma(\varepsilon)$. Denote by $\tilde \gamma$ the $g_F$-geodesic segment connecting $\gamma(-\varepsilon)$ and $\gamma(\varepsilon)$, and denote by $C$ its  $g_F$-midpoint, see Fig. 1.

\noindent\begin{minipage}{.49\textwidth}
\includegraphics[width=\textwidth]{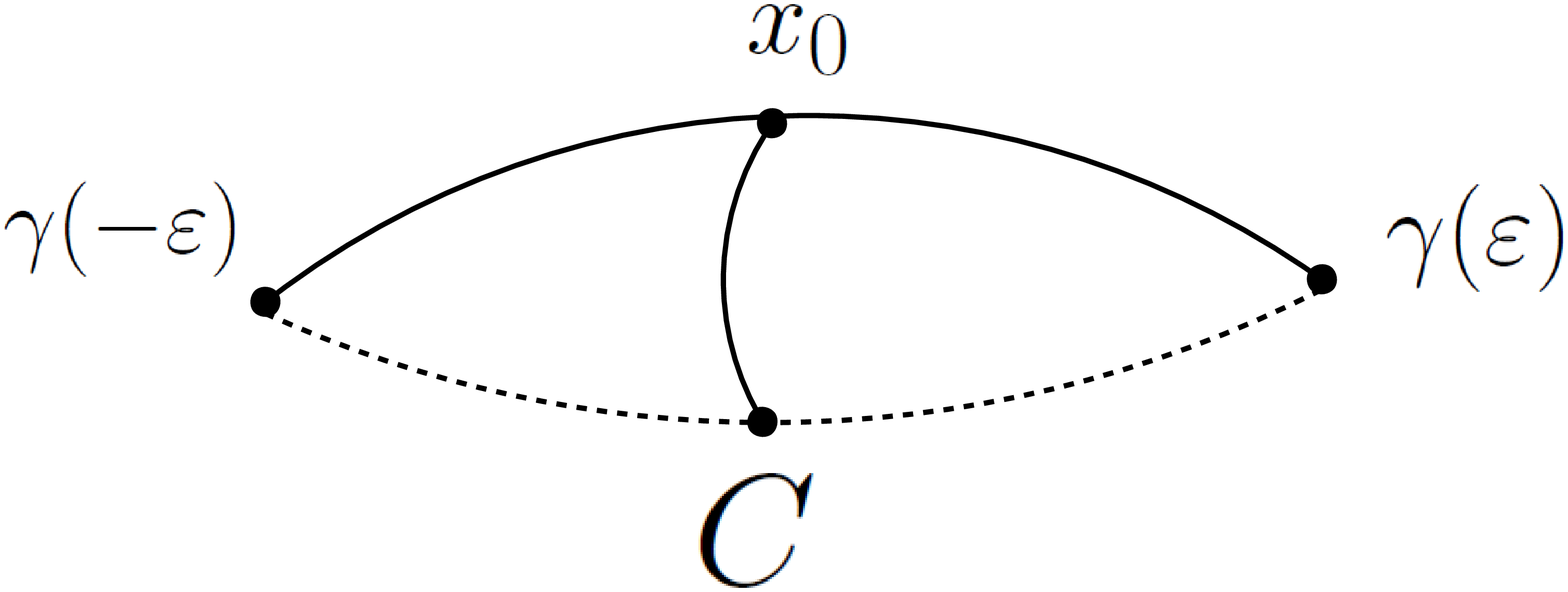}   \end{minipage}\begin{minipage}{.4\textwidth}
%\begin{figure} 
{\bf Fig.1.} {If the point $C$ is different from $x_0$, then the geodesic connecting the point $x_0$ with the point $C$ is preserved by the  isometry $\phi^{U_i}$ which contradicts that its initial velocity vector is not preserved by $d\phi_{x_0}^{U_i}$  }\label{pic5} %\end{figure}
 \end{minipage}

\vspace{1ex} 
Let us show that $C=x_0$. Assume it is not the case.  Denote by $Z\in T_{x_0}M$ the initial velocity vector of the $F$-geodesic connecting $x_0$ and $C$. Now, take  the isometry $\phi^{U_i}$ such that $d\phi^{U_i}_{x_0}(Z)\ne Z$. 
The isometry is well-defined in the r-ball around $x_0$, since all  isometries $\phi^{U_1},...,\phi^{U_N}$ are well-defined there. 
It sends  the geodesic segment $\gamma$ to itself with reversing the time, since $d\phi^{U_i}(\dot\gamma(0))=  -\dot \gamma(0)$. It therefore sends $\gamma(\varepsilon)$ to $\gamma(-\varepsilon)$ and $\gamma(-\varepsilon)$ to $\gamma(\varepsilon)$. Then, it sends the segment $\tilde \gamma$ to itself, which implies it sends its midpoint, which is $C$,  to itself. Then, it sends the geodesic connecting $x_0$ and $C$ to itself, which implies that its differential $d\phi_{x_0}^{U_i}$
 sends $Z$  to itself, which contradicts the choice of $\phi^{U_i}$.

Thus, the geodesic $\tilde \gamma$ passes through $x_0$ which implies that it is a straight line in our coordinate system. By construction it has three intersections with $\gamma$, namely $\gamma(-\varepsilon)$, $x_0$  and $\gamma(\varepsilon)$, as we want. 
Finally, $G^i(x_0, V)=0$.

Thus, the metric $F$  is Berwaldian and the associated connection is the Levi-Civita connection of $g_F$.   Since the associated connection is the Levi-Civita connection of a locally symmetric Riemannian manifold, the parallel transport in this connection commutes with the geodesic reflection. Then, the geodesic reflection is an isometry for $F$ which implies that $F$ is locally symmetric.
 Theorem \ref{thm2} is proved.

\noindent Shaoqiang Deng:  School of Mathematical Sciences and LPMC,
                            Nankai University,
                            Tianjin 300071, China.  Email: \url{dengsq@nankai.edu.cn}.  \\   Vladimir S. Matveev: Institut f\"ur  Mathematik, 
Friedrich-Schiller Universit\"at Jena,    
07737    
Germany. Email: \url{vladimir.matveev@uni-jena.de}.

\end{document}